\documentclass{amsart}

\usepackage{amssymb}
\usepackage{needspace}
\usepackage{dsfont}
\usepackage{xcolor}
\usepackage{graphicx}
\usepackage{hyperref}
\hypersetup{colorlinks, linkcolor={red!50!black}, citecolor={blue!50!black},
  urlcolor={blue!80!black}}
\usepackage{doi}

\usepackage[T1]{fontenc}
\usepackage{booktabs}
\usepackage{tabularx}
\usepackage{array}
\usepackage{listings}
\usepackage{seqsplit}

\newcommand{\pkg}{\texttt{quantum-group}}
\newcommand{\code}[1]{\texttt{\begingroup\def\_{\textunderscore\allowbreak}#1\endgroup}}
\newcommand{\Uq}{U_q(\mathfrak{sl}_2)}
\newcommand{\GL}{GL_q(2|1)}

\lstdefinestyle{output}{frame=none, basicstyle=\ttfamily\footnotesize\color{black!80},
  aboveskip=0.2em, belowskip=0.6em}

\usepackage[labelfont=sc, labelsep=period, font=small]{caption}
\makeatletter
\def\@defaultbiblabelstyle#1{[#1]}
\makeatother

\makeatletter
\renewcommand\subsection{\@startsection{subsection}{2}{\z@}%
  {.7\linespacing\@plus\linespacing}{.3\linespacing}{\normalfont\bfseries}}
\makeatother

\DeclareRobustCommand{\qgname}{\texttt{quantum-group}}

\hypersetup{pdftitle={quantum-group: Exact symbolic verification of R-matrix and Yang-Baxter identities for Uq(sl2) and GLq(2|1)},
  pdfauthor={Taha Berk Terekli}}

\begin{document}

\title[\qgname: exact verification of R-matrix identities]{\qgname: Exact symbolic verification of R-matrix and Yang--Baxter identities for $U_q(\mathfrak{sl}_2)$ and $GL_q(2|1)$}

\author{Taha Berk Terekli}
\address{Department of Mathematics, Y{\i}ld{\i}z Technical University, Istanbul, T\"urkiye}
\email{terekli@tahaberk.com}
\urladdr{https://orcid.org/0009-0004-8266-1116}


\begin{abstract}
Explicit R-matrices of quantum groups are routinely transcribed from the literature, yet whether they satisfy an identity depends on basis order, tensor placement, the coproduct and, for superalgebras, the Koszul sign rule. \pkg{} is a Python package that builds such matrices exactly with SymPy and returns every identity as an exact residual matrix, whose nonzero entries locate the convention that disagrees. It covers $\Uq$ modules, tensor products and the fundamental R-matrix, and the graded Yang--Baxter equation for $\GL$. Applied to the thesis it was written for, it exposed an R-matrix that satisfies the Yang--Baxter equation but intertwines the opposite coproduct.
\end{abstract}

\maketitle

\section{Introduction}
\label{sec:intro}

Quantum groups are one-parameter deformations of classical symmetry algebras. Their representations yield solutions of the quantum Yang--Baxter equation (QYBE), which is why they appear in exactly solvable models of statistical mechanics, in knot invariants and in representation theory~\cite{drinfeld1987,jimbo1985,jimbo1986,kassel1995}. Much of the working knowledge in this area is recorded as explicit matrices whose entries are rational functions of the deformation parameter $q$: an R-matrix on $V\otimes V$, the generators of a module, a change of basis to irreducible summands. Before such a matrix can be reused it has to be checked, and checking it means multiplying $9\times9$, $27\times27$ or $81\times81$ matrices whose meaning depends on conventions that texts state in different ways, or not at all.

The characteristic mistakes are not arithmetic ones. A matrix is transcribed in a different basis order from the one its source used. $R_{13}$ is formed by conjugating with the wrong permutation. An R-matrix is paired with the opposite of the coproduct it actually intertwines; this still satisfies the QYBE, so the most common test does not detect it. In a $\mathbb{Z}_2$-graded setting an ordinary flip is used where the Koszul sign rule requires a super-permutation. And a check that is run at a floating-point value of $q$ can hide cancellations or produce spurious nonzero entries. Each of these errors produces output that looks plausible.

\pkg{}~\cite{quantumgroup111} addresses this with a deliberately narrow design. Every object is an ordinary exact \code{sympy.Matrix}~\cite{meurer2017sympy}; conventions (basis order, coproduct, parities) are fixed in documented functions; and every identity is exposed as a \emph{residual}, the exact matrix $\mathrm{LHS}-\mathrm{RHS}$, which is then tested entry by entry. For $\Uq$ the package constructs the irreducible modules $V_n$, checks the defining relations and generator-level Hopf identities, builds tensor products through the coproduct, computes Clebsch--Gordan highest-weight vectors, and verifies the fundamental R-matrix: QYBE, braid and Hecke relations and coproduct intertwining. For the quantum supergroup $\GL$ introduced by Çelik and Çelik~\cite{celik2021}, whose third basis vector is odd, it constructs the published $9\times9$ R-matrix and verifies the graded Yang--Baxter equation on three factors and on every triple of factors in $V^{\otimes4}$.

The software was developed for and used in the author's undergraduate thesis~\cite{terekli2026thesis}, where it produced the thesis's tables, matrix figures and Yang--Baxter verifications. After the thesis was submitted in June 2026, the software was used to audit it: an explicit intertwining residual showed that the thesis's R-matrix intertwines the opposite of the coproduct used in the text, an error that passing QYBE, braid and Hecke checks had not revealed (Sec.~\ref{sec:ex-coproduct}). The correction is recorded in the repository (\code{thesis/ERRATA.md}, \code{AUDIT.md}), and the check is now part of the API and the test suite. No new mathematical results are claimed; the contribution is the software and its verification method.

The intended users are researchers who need to confirm a matrix identity from the literature in their own conventions: mathematical physicists working with R-matrices and integrable models, representation theorists and computer-algebra users who want an inspectable reference calculation, and students reproducing standard quantum-group computations.

\par\needspace{3\baselineskip}\smallskip\noindent\emph{Related software.} GAP's QuaGroup package~\cite{gap,quagroup,degraaf2001} and SageMath~\cite{sagemath}, which provides a \code{QuantumGroup} interface to it, work with the quantized enveloping algebra $U_q(\mathfrak{g})$ of any semisimple Lie algebra $\mathfrak{g}$: they construct PBW-type and canonical bases, highest-weight modules and their R-matrices; SageMath also implements crystals for the quantum superalgebras $U_q(\mathfrak{gl}(m|n))$~\cite{bkk2000}. These systems are far broader, and for general $U_q(\mathfrak{g})$ work they are the appropriate tools. \pkg{} does not replace them. Its added value is that of a small verification layer for concrete matrix calculations: it keeps the basis order and coproduct of a particular text, returns the matrices and residuals themselves rather than only a verdict, runs in an ordinary Python environment, encodes the graded sign rule once for every tensor placement, and includes a graded $\GL$ example that lies outside QuaGroup's documented $U_q(\mathfrak{g})$ scope. The two are complementary: Sec.~\ref{sec:crosscheck} uses \pkg{} to test matrices computed by QuaGroup and to identify the conventions in which the systems differ.

\par\needspace{3\baselineskip}\smallskip\noindent\emph{Outline.} Section~\ref{sec:background} collects the definitions the package implements, and Sec.~\ref{sec:design} describes how it works and why it is built this way. Section~\ref{sec:usage} explains installation and the public interface, Sec.~\ref{sec:examples} gives worked examples, including the thesis audit, and Sec.~\ref{sec:tests} documents testing, a cross-check against GAP's QuaGroup and performance. Section~\ref{sec:scope} states the scope and the ways the code can be reused or extended. Appendix~\ref{app:api} is a reference of the public interface.

\section{Mathematical background}
\label{sec:background}

This section collects the definitions that the package implements, in the conventions it uses. Standard references are Refs.~\cite{drinfeld1987,jimbo1985,jimbo1986,kassel1995}; the graded case follows Ref.~\cite{celik2021}.

\subsection{The algebra \texorpdfstring{$\Uq$}{Uq(sl2)} and its modules}
\label{sec:uq}

Let $q$ be an indeterminate, or a nonzero complex number that is not a root of unity. $\Uq$ is the algebra generated by $E$, $F$, $K$ and $K^{-1}$ with the relations
\begin{equation}
\begin{gathered}
\text{(R1)}\ \ KK^{-1}=K^{-1}K=1,\qquad
\text{(R2)}\ \ KEK^{-1}=q^{2}E,\\
\text{(R3)}\ \ KFK^{-1}=q^{-2}F,\qquad
\text{(R4)}\ \ EF-FE=\frac{K-K^{-1}}{q-q^{-1}} .
\end{gathered}
\label{eq:relations}
\end{equation}
It is a Hopf algebra with coproduct, counit and antipode
\begin{gather}
\begin{gathered}
\Delta(E)=E\otimes 1+K\otimes E,\qquad \Delta(F)=F\otimes K^{-1}+1\otimes F,\\
\Delta(K^{\pm1})=K^{\pm1}\otimes K^{\pm1},
\end{gathered}
\label{eq:coproduct}\\
\begin{gathered}
\varepsilon(E)=\varepsilon(F)=0,\qquad \varepsilon(K^{\pm1})=1,\\
S(E)=-K^{-1}E,\qquad S(F)=-FK,\qquad S(K^{\pm1})=K^{\mp1}.
\end{gathered}
\label{eq:hopf}
\end{gather}
The opposite coproduct is $\Delta^{\mathrm{op}}=P\circ\Delta$, where $P$ exchanges the two tensor factors; on a tensor product of modules, $\Delta^{\mathrm{op}}(X)=P\,\Delta(X)\,P$. With the $q$-integers
\begin{equation}
[k]=\frac{q^{k}-q^{-k}}{q-q^{-1}}=q^{k-1}+q^{k-3}+\dots+q^{1-k},\qquad [k]!=[1][2]\cdots[k],
\end{equation}
the irreducible module $V_n$ of dimension $n+1$ has a basis $v_0,\dots,v_n$ on which
\begin{equation}
\begin{gathered}
K v_k=q^{n-2k}v_k,\qquad F v_k=v_{k+1}\ \ (F v_n=0),\\
E v_k=[k]\,[n-k+1]\,v_{k-1}\ \ (E v_0=0).
\end{gathered}
\label{eq:module}
\end{equation}
A tensor product $V\otimes W$ of modules is a module through $\Delta$. For generic $q$ the $q$-analogue of the Clebsch--Gordan rule holds, $V_m\otimes V_n\cong\bigoplus_{k}V_k$ with $k=|m-n|,|m-n|+2,\dots,m+n$, and each summand is generated by a \emph{highest-weight vector}, a $K$-eigenvector annihilated by $E$.

\subsection{R-matrices and Yang--Baxter-type identities}
\label{sec:ybe}

Let $V$ have dimension $d$ and let $R$ be an operator on $V\otimes V$, i.e.\ a $d^2\times d^2$ matrix. On $V^{\otimes3}$ one writes $R_{12}=R\otimes 1$, $R_{23}=1\otimes R$ and $R_{13}=P_{23}R_{12}P_{23}$, where $P_{23}$ exchanges the second and third factors. The quantum Yang--Baxter equation and, for the braided form $\check R=PR$, the braid relation read
\begin{equation}
R_{12}R_{13}R_{23}=R_{23}R_{13}R_{12},\qquad
\check R_{12}\check R_{23}\check R_{12}=\check R_{23}\check R_{12}\check R_{23};
\label{eq:qybe}
\end{equation}
the two are equivalent~\cite{kassel1995}. For the fundamental representation, $\check R$ moreover satisfies the Hecke relation
\begin{equation}
\check R-\check R^{-1}=(q-q^{-1})\,\mathds{1},\qquad\text{equivalently}\qquad (\check R-q)(\check R+q^{-1})=0 .
\label{eq:hecke}
\end{equation}
Neither relation refers to the Hopf structure. What ties an R-matrix on $V\otimes V$ to a coproduct is the intertwining condition
\begin{equation}
R\,\Delta(X)=\Delta^{\mathrm{op}}(X)\,R\qquad\text{for all generators }X,
\label{eq:intertwining}
\end{equation}
which implies that $\check R=PR$ commutes with the action of $\Delta(X)$ and is therefore a module endomorphism of $V\otimes V$. A matrix can satisfy Eqs.~\eqref{eq:qybe} and~\eqref{eq:hecke} and yet fulfil Eq.~\eqref{eq:intertwining} only for $\Delta^{\mathrm{op}}$ in place of $\Delta$; Section~\ref{sec:ex-coproduct} shows such a case.

\subsection{\texorpdfstring{$\mathbb{Z}_2$}{Z2}-graded tensor products and \texorpdfstring{$\GL$}{GLq(2|1)}}
\label{sec:grading}

A super vector space $V=V_{\bar0}\oplus V_{\bar1}$ comes with a basis $e_i$ of homogeneous vectors of parity $p(i)\in\{0,1\}$; the basis vector $e_a\otimes e_b$ of $V\otimes V$ has parity $p(a)+p(b)\bmod 2$. An operator is \emph{even} if it maps even vectors to even ones and odd to odd ones, i.e.\ if its matrix element between basis vectors of different parity vanishes. In tensor products the Koszul sign rule applies: exchanging homogeneous elements of parities $a$ and $b$ produces the sign $(-1)^{ab}$. Accordingly the flip $P$ is replaced by the super-permutation
\begin{equation}
P_s(e_i\otimes e_j)=(-1)^{p(i)p(j)}\,e_j\otimes e_i ,
\label{eq:superperm}
\end{equation}
and an even operator on $V\otimes V$ is placed on factors $(i,j)$ of $V^{\otimes n}$ by conjugation with graded rather than ordinary flips. With $R_{13}$ formed in this way, the \emph{graded} Yang--Baxter equation has the form of Eq.~\eqref{eq:qybe}.

For $\GL$, $V=\mathbb{C}^{2|1}$ with basis $e_0,e_1$ (even) and $e_2$ (odd); Ref.~\cite{celik2021} numbers these vectors from 1. In the basis $e_a\otimes e_b$, ordered by $3a+b$, the R-matrix printed on page~261 of Ref.~\cite{celik2021}, as implemented in \code{R\_matrix\_GLq21}, is
\begin{equation}
R=\begin{pmatrix}
1&0&0&0&0&0&0&0&0\\
0&q^2&0&1-q^2&0&0&0&0&0\\
0&0&q&0&0&0&1-q^2&0&0\\
0&0&0&1&0&0&0&0&0\\
0&0&0&0&1&0&0&0&0\\
0&0&0&0&0&q&0&1-q^2&0\\
0&0&0&0&0&0&q&0&0\\
0&0&0&0&0&0&0&q&0\\
0&0&0&0&0&0&0&0&q^2
\end{pmatrix}.
\label{eq:rgl21}
\end{equation}
Its three off-diagonal entries couple $e_0\otimes e_1$ with $e_1\otimes e_0$, $e_0\otimes e_2$ with $e_2\otimes e_0$, and $e_1\otimes e_2$ with $e_2\otimes e_1$, pairs of equal parity, so $R$ is even. It satisfies the graded Yang--Baxter equation but not the ungraded one (Sec.~\ref{sec:ex-first}).

\section{Design and workings}
\label{sec:design}

\begin{figure}[t]
  \centering
  \includegraphics[width=\textwidth]{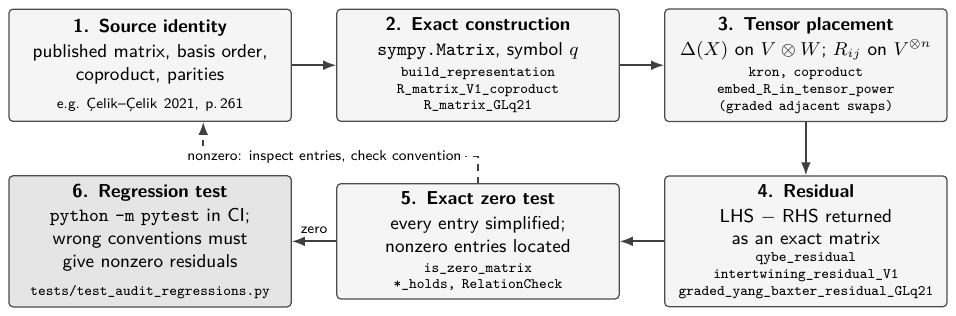}
  \caption{Verification workflow in \pkg. A published identity and its conventions (1) are turned into exact SymPy matrices (2) and placed on tensor products through the coproduct or through graded adjacent swaps (3). The difference of the two sides is returned as an exact residual matrix (4) and tested entry by entry (5); a nonzero residual points to the entries, and hence the convention, that disagree. Zero residuals, and nonzero residuals for deliberately wrong conventions, are kept as regression tests (6). Function names are those of the package; \code{kron} is imported from \code{quantum\_group.linalg}, all others from the package namespace.}
  \label{fig:workflow}
\end{figure}

\pkg{} is a pure-Python package of fifteen modules (about 2{,}100 lines, including docstrings) with no compiled extensions. Figure~\ref{fig:workflow} shows the verification workflow that the modules implement; the design decisions behind each step are described below.

\subsection{Explicit matrices rather than an abstract algebra engine}

The package contains no noncommutative rewriting system; the symbols in \code{generators.py} only display the relations. Every check runs on concrete representation matrices. This keeps each object directly comparable, entry by entry, with a printed matrix, and makes every intermediate result an ordinary SymPy object that users can print, substitute into or pass to other SymPy code. The price is that a verification holds for the representations on which it is run, and that the cost grows with the dimension of the tensor product (Sec.~\ref{sec:benchmarks}). Identities that need a genuinely noncommutative graded algebra, such as the Gauss decomposition of $\GL$, are omitted rather than approximated with commuting symbols.

\subsection{Residual-first verification}
\label{sec:residual}

Each check follows one pattern: build the left-hand side, build the right-hand side, return their difference, and decide whether every entry is zero. Functions named \code{*\_residual} (\code{qybe\_residual}, \code{braid\_relation\_residual}, \code{intertwining\_residual\_V1}, \code{graded\_yang\_baxter\_residual\_GLq21}, and others) return the exact difference matrix; the corresponding \code{*\_holds} functions apply \code{is\_zero\_matrix} to it; and \code{verify\_on\_representation} returns, for each defining relation, a \code{RelationCheck} record holding both the Boolean and the residual. Returning the residual rather than only a Boolean is what makes a failed check informative: the positions and values of the nonzero entries identify which basis vectors and which convention disagree (Sec.~\ref{sec:examples}).

The zero test simplifies each entry with \code{sympy.simplify}. A \code{True} result means that every entry was reduced to the exact value 0; a \code{False} result is not by itself a proof of inequality, because \code{simplify} is not a decision procedure. For the rational functions of $q$ produced by the package, a nonzero entry can be confirmed with \code{sympy.cancel}, which the regression tests do.

\subsection{Conventions for \texorpdfstring{$\Uq$}{Uq(sl2)}}
\label{sec:conventions}

For a highest weight $n$, the function \code{build\_representation(n)} returns the $(n+1)$-dimensional module $V_n$ as a \code{Representation} record with explicit matrices $E$, $F$, $K$ and $K^{-1}$ in the basis $v_0,\dots,v_n$, as in Eq.~\eqref{eq:module}: $K$ is diagonal, $F$ is a lower shift and $E$ carries the products of $q$-integers. \code{coproduct} and \code{tensor\_product} realize the coproduct~\eqref{eq:coproduct} with Kronecker products in the lexicographic basis order (left factor slow, right factor fast). \code{find\_highest\_weight\_vectors} solves for the vectors annihilated by $E$ in each weight space of $V_m\otimes V_n$, giving one highest-weight vector per Clebsch--Gordan summand.

\par\needspace{3\baselineskip}\smallskip\noindent\emph{Two R-matrix conventions.} On $V_1\otimes V_1$, in the basis $v_0\otimes v_0,\,v_0\otimes v_1,\,v_1\otimes v_0,\,v_1\otimes v_1$, the package provides
\begin{equation}
R=\begin{pmatrix} q&0&0&0\\ 0&1&q-q^{-1}&0\\ 0&0&1&0\\ 0&0&0&q\end{pmatrix},
\qquad
R_{21}=\begin{pmatrix} q&0&0&0\\ 0&1&0&0\\ 0&q-q^{-1}&1&0\\ 0&0&0&q\end{pmatrix},
\label{eq:rmatrices}
\end{equation}
where $R_{21}=PRP$ and $P$ is the flip of the two factors. \code{R\_matrix\_V1} returns the historical upper-triangular matrix $R$, the image, rescaled by $q^{1/2}$, of Drinfeld's universal R-matrix~\cite{drinfeld1987} for the coproduct opposite to Eq.~\eqref{eq:coproduct}; accordingly it satisfies $R\,\Delta^{\mathrm{op}}(X)=\Delta(X)\,R$, i.e.\ it intertwines the \emph{opposite} of the coproduct~\eqref{eq:coproduct}. Rather than change its output, which would silently alter results already computed with it, version 1.1.0 kept it unchanged and added \code{R\_matrix\_V1\_coproduct} ($R_{21}$), \code{R\_check\_V1\_coproduct} ($PR_{21}$) and \code{intertwining\_residual\_V1}, which satisfy $R_{21}\Delta(X)=\Delta^{\mathrm{op}}(X)R_{21}$. Both conventions satisfy the QYBE, and both braided forms satisfy the braid relation and the Hecke relation $\check R-\check R^{-1}=(q-q^{-1})\,\mathds{1}$, with eigenvalues $q$ (multiplicity three) and $-q^{-1}$; they differ in which tensor-product submodules are the braid eigenspaces. Keeping both, with the convention stated in each docstring, lets earlier calculations be reproduced exactly while new work uses the compatible matrices. The same policy applies to names: a renamed function keeps its old name as a thin wrapper that emits a \code{DeprecationWarning}.

\subsection{Graded tensor embeddings for \texorpdfstring{$\GL$}{GLq(2|1)}}
\label{sec:graded}

For $\GL$ the basis parities $(0,0,1)$ are an explicit argument, and the graded flip $P_s$ of Eq.~\eqref{eq:superperm} is built once, by \code{super\_permutation\_matrix}. \code{embed\_R\_in\_tensor\_power} places a $d^2\times d^2$ operator on factors $(i,j)$ of $V^{\otimes n}$ by first embedding it on the adjacent pair $(i,i+1)$ with Kronecker products and then conjugating by graded adjacent swaps $S_k=\mathds{1}^{\otimes k}\otimes P_s\otimes \mathds{1}^{\otimes(n-k-2)}$ until the second factor reaches position $j$. The Koszul sign rule is therefore encoded in a single matrix and reused for every placement, instead of being written out separately for $R_{13}$, $R_{14}$, $R_{24}$ and so on. The construction is correct only for even (parity-preserving) operators, and the function checks this before embedding. Because the parity list is an argument, the same function also embeds operators on other super vector spaces, and an all-even parity list reproduces the ungraded placement. \code{R\_matrix\_GLq21} returns the $9\times9$ R-matrix~\eqref{eq:rgl21} exactly as printed on page~261 of Ref.~\cite{celik2021}, and \code{graded\_yang\_baxter\_residual\_GLq21} forms $R_{12}R_{13}R_{23}-R_{23}R_{13}R_{12}$ on $V^{\otimes3}$.

\subsection{Exact parameters}

The default parameter is a nonzero SymPy symbol, and the identities are rational in $q$, holding wherever their denominators are defined. Integers and exact rationals passed as $q$ are converted to exact SymPy numbers, so \code{build\_representation(2, 2)} yields entries such as $5/2$ rather than $2.5$; only floating-point values supplied by the user are approximate. $q$-integers and $q$-binomials are evaluated as Laurent polynomials, not as quotients, so that specializations with removable singularities, such as $q=-1$ or roots of unity, give the correct values instead of $0/0$.

\subsection{Validation and failure modes}

Public functions reject invalid input with specific errors rather than returning a misleading matrix: non-integral or negative highest weights, sizes and indices raise \code{TypeError} or \code{ValueError}; \code{verify\_on\_representation} rejects non-square or mismatched generator matrices and the singular values $q\in\{0,1,-1\}$ of the defining commutator quotient; \code{qybe\_residual} rejects matrices that are not $d^2\times d^2$; and \code{embed\_R\_in\_tensor\_power} rejects invalid positions, parity lists and odd operators. Symbolic input is not restricted: any SymPy expression can be used as $q$.

\subsection{Package structure}

The modules form three layers. A small foundation provides the default symbol $q$ (declared nonzero), $q$-integers, $q$-factorials and $q$-binomials (\code{utils}); the Kronecker product \code{kron}/\code{kron\_list} and the zero-residual test \code{is\_zero\_matrix} (\code{linalg}); and internal argument checks (\code{\_validation}). Construction modules build the matrices: $\Uq$ modules (\code{representations}), coproduct, counit and antipode (\code{hopf}), tensor products and Clebsch--Gordan vectors (\code{tensor}), $\Uq$ R-matrices (\code{r\_matrix}) and the $\GL$ R-matrix, super-permutation and embeddings (\code{supergroup\_gl21}). Verification functions sit next to the constructions they check and all use the same foundation. Auxiliary modules provide classical and root-of-unity specializations (\code{limits}), the crystal graph $B(n)$ as a NetworkX~\cite{hagberg2008networkx} graph (\code{crystal}), Matplotlib~\cite{hunter2007matplotlib} weight and crystal plots (\code{visualization}) and a convenience facade, \code{QuantumGroupSL2}. The main functions are re-exported from the package namespace and listed in \code{\_\_all\_\_}.

\section{Installation and usage}
\label{sec:usage}

\subsection{Requirements, installation and availability}
\label{sec:install}

\pkg{} requires CPython 3.10 or later and runs on any operating system; it has no compiled code and needs no special hardware. It is tested in continuous integration on Linux, macOS and Windows (Sec.~\ref{sec:tests}) and locally on Ubuntu 24.04 (x86-64). The runtime dependencies, installed automatically, are SymPy $\geq$ 1.10~\cite{meurer2017sympy}, NetworkX $\geq$ 2.6~\cite{hagberg2008networkx} and Matplotlib $\geq$ 3.5~\cite{hunter2007matplotlib}; continuous integration exercises the oldest supported versions (SymPy 1.10.1, NetworkX 2.6.3, Matplotlib 3.5.3). pytest~\cite{pytest} is needed only to run the test suite (extra \code{[test]}).

The version described here is 1.1.1, tagged as the GitHub Release \code{v1.1.1} (commit \code{7cb3877}). Its package code and test suite are those of version 1.1.0; release 1.1.1 adds the listings of this paper with their expected outputs, the data of the cross-check in Sec.~\ref{sec:crosscheck}, and continuous-integration jobs that run them. It is installed from the Python Package Index, or directly from the tag:
\begin{lstlisting}[language=bash, frame=none]
python -m pip install quantum-group==1.1.1
python -m pip install \
    "git+https://github.com/TerekliTahaBerk/quantum-groups@v1.1.1"
python -c "import quantum_group; print(quantum_group.__version__)"   # 1.1.1
\end{lstlisting}
To run the tests and benchmarks, clone the repository and install it with the test extra:
\begin{lstlisting}[language=bash, frame=none]
git clone https://github.com/TerekliTahaBerk/quantum-groups.git
cd quantum-groups && python -m pip install -e ".[test]"
python -m pytest                                          # 202 tests
python benchmarks/benchmark.py --repeats 3 --max-n 5      # Table 3
\end{lstlisting}

\par\needspace{3\baselineskip}\smallskip\noindent\emph{Code availability.} Repository: \url{https://github.com/TerekliTahaBerk/quantum-groups} (MIT licence). Release: \url{https://github.com/TerekliTahaBerk/quantum-groups/releases/tag/v1.1.1}. Package: \url{https://pypi.org/project/quantum-group/1.1.1/}. Archive: Zenodo, \doi{10.5281/zenodo.23039671}~\cite{quantumgroup111}, a source archive of the GitHub Release; version 1.1.0 is archived as \doi{10.5281/zenodo.23015576}~\cite{quantumgroup110}. Listing~\ref{lst:example} is \code{examples/sample\_verification.py}, and Listings~\ref{lst:coproduct}--\ref{lst:crosscheck} are in \code{scipost/examples/} and \code{scipost/crosscheck/}; all six are part of the release, each with its expected output.

\subsection{Overview of the public interface}

Table~\ref{tab:api} lists the main functions by task; Appendix~\ref{app:api} gives their arguments and return values. Functions that construct $q$-dependent objects accept an optional parameter \code{q\_sym} (default: the package symbol \code{q}); all return ordinary SymPy objects or plain Python records of them, and each states its convention in the docstring.

\begin{table}[t]
\centering
\caption{Main public functions of \pkg{} 1.1.1, grouped by task. All are importable from the package namespace.}
\label{tab:api}
\small
\begin{tabularx}{\textwidth}{@{}>{\raggedright\arraybackslash}p{0.25\textwidth}>{\raggedright\arraybackslash}X@{}}
\toprule
Task & Functions \\
\midrule
$\Uq$ modules and relations & \texttt{build\_representation}, \texttt{verify\_on\_representation}, \texttt{all\_relations\_hold}, \texttt{symbolic\_relations} \\
Hopf structure & \texttt{coproduct}, \texttt{counit}, \texttt{antipode}, \texttt{verify\_coassociativity}, \texttt{verify\_antipode}, \texttt{verify\_all\_hopf\_axioms} \\
Tensor products & \texttt{tensor\_product}, \texttt{find\_highest\_weight\_vectors}, \texttt{cg\_summands}, \texttt{cg\_decomposition\_summary}, \texttt{swap\_matrix}, \texttt{kron\_list} \\
$\Uq$ R-matrices & \texttt{R\_matrix\_V1}, \texttt{R\_matrix\_V1\_coproduct}, \texttt{R\_check\_V1}, \texttt{R\_check\_V1\_coproduct}, \texttt{R\_check\_eigenvalues} \\
Generic residuals & \texttt{qybe\_residual}, \texttt{qybe\_holds}, \texttt{braid\_relation\_residual}, \texttt{braid\_relation\_holds}, \texttt{hecke\_skein\_relation\_check}, \texttt{intertwining\_residual\_V1}, \texttt{intertwining\_holds\_V1}, \texttt{is\_zero\_matrix} \\
$\GL$ and grading & \texttt{R\_matrix\_GLq21}, \texttt{super\_parity\_gl21}, \texttt{super\_permutation\_matrix}, \texttt{embed\_R\_in\_tensor\_power}, \texttt{all\_Rij\_GLq21}, \texttt{graded\_yang\_baxter\_residual\_GLq21}, \texttt{graded\_yang\_baxter\_holds\_GLq21}, \texttt{local\_ybe\_on\_four\_tensor\_GLq21}, \texttt{braid\_far\_commutativity\_residual\_GLq21} \\
Specializations & \texttt{q\_integer}, \texttt{q\_factorial}, \texttt{q\_binomial}, \texttt{classical\_limit}, \texttt{root\_of\_unity\_substitution} \\
Crystals and plots & \texttt{build\_crystal}, \texttt{e\_tilde}, \texttt{f\_tilde}, \texttt{plot\_weight\_diagram}, \texttt{plot\_crystal\_graph} \\
\bottomrule
\end{tabularx}
\end{table}

\subsection{Documentation and support}

Documentation consists of the README (installation, quickstart, conventions, limitations), module and function docstrings that state the conventions used (basis order, coproduct, parities) and the relation each check tests, \code{CONTRIBUTING.md}, \code{CHANGELOG.md}, an exploratory Jupyter notebook, and \code{MANUSCRIPT\_CODE\_MAPPING.md}, which links the computational claims of the thesis to their implementation and tests. Questions, bug reports and feature requests are handled through the GitHub issue tracker (\url{https://github.com/TerekliTahaBerk/quantum-groups/issues}). \code{CONTRIBUTING.md} describes the development set-up, the validation command and the conventions for new checks (a \code{*\_residual} function with a \code{*\_holds} companion, input validation with a test for each error path, and a stated reference for any convention). The package is maintained by its author on a best-effort basis, with no institutional support or guaranteed response time; changes are recorded in \code{CHANGELOG.md} and versions follow semantic versioning.

\section{Worked examples}
\label{sec:examples}

The outputs below were produced with \pkg{} 1.1.1, SymPy 1.14 and Python 3.11; each is stored next to its script and is compared automatically with a fresh run.

\subsection{A first verification session}
\label{sec:ex-first}

Listing~\ref{lst:example} (\code{examples/sample\_verification.py}) checks an installation in a few seconds. It verifies the defining relations on $V_2$, the QYBE and coproduct intertwining of $R_{21}$, and the graded Yang--Baxter equation of $\GL$; it then repeats the last check with the $\GL$ R-matrix placed by ordinary, all-even swaps, as a negative control.

\begin{lstlisting}[language=Python, caption={First verification session: \code{examples/sample\_verification.py}.}, label={lst:example}]
from quantum_group import (
    build_representation, verify_on_representation, R_matrix_V1_coproduct,
    qybe_holds, intertwining_holds_V1, R_matrix_GLq21,
    embed_R_in_tensor_power, graded_yang_baxter_residual_GLq21,
    is_zero_matrix)

rep = build_representation(2)            # V_2 of U_q(sl_2): 3x3 matrices
print("E on V_2:", rep.E.tolist())
checks = verify_on_representation(rep.E, rep.F, rep.K, rep.K_inv)
print("relations:", {k: c.holds for k, c in checks.items()})

R = R_matrix_V1_coproduct()              # R_21, matches the package coproduct
print("QYBE:", qybe_holds(R))
print("intertwines Delta:",
      [intertwining_holds_V1(X) for X in ("E", "F", "K", "K_inv")])

Y = graded_yang_baxter_residual_GLq21()  # GL_q(2|1), parities (0, 0, 1)
print("graded YBE residual:", Y.shape, "zero:", is_zero_matrix(Y))

# Negative control: the same R placed with ordinary (all-even) swaps.
R12, R13, R23 = (embed_R_in_tensor_power(R_matrix_GLq21(), ij, 3, [0, 0, 0])
                 for ij in [(0, 1), (0, 2), (1, 2)])
bad = (R12 * R13 * R23 - R23 * R13 * R12).expand()
print("ungraded residual, nonzero entries:",
      {(i, j): x.factor()
       for (i, j), x in sorted(bad.todok().items()) if x != 0})
\end{lstlisting}
Output (\code{examples/sample\_verification\_expected.txt}; the last line is wrapped here):
\begin{lstlisting}[style=output]
E on V_2: [[0, q + 1/q, 0], [0, 0, q + 1/q], [0, 0, 0]]
relations: {'R1': True, 'R2': True, 'R3': True, 'R4': True}
QYBE: True
intertwines Delta: [True, True, True, True]
graded YBE residual: (27, 27) zero: True
ungraded residual, nonzero entries: {(8, 24): 2*q**2*(q - 1)**2*(q + 1)**2,
  (17, 25): 2*q**2*(q - 1)**2*(q + 1)**2}
\end{lstlisting}
The negative control shows why the grading matters: with ordinary swaps, exactly two entries of the $27\times27$ residual are nonzero, both equal to $2q^2(q-1)^2(q+1)^2$. Rows and columns are indexed lexicographically, $e_a\otimes e_b\otimes e_c\mapsto 9a+3b+c$, so entry $(8,24)$ couples $e_0\otimes e_2\otimes e_2$ with $e_2\otimes e_2\otimes e_0$; both involve the odd vector $e_2$ twice, which is exactly where the Koszul sign $(-1)^{p(i)p(j)}$ differs from $+1$.

\subsection{Application: which coproduct does an R-matrix intertwine?}
\label{sec:ex-coproduct}

This is the check that exposed the convention error in the author's thesis~\cite{terekli2026thesis}. The thesis used the coproduct~\eqref{eq:coproduct} together with the upper-triangular matrix $R$ of Eq.~\eqref{eq:rmatrices}. $R$ satisfies the QYBE and its braided form satisfies the braid and Hecke relations, so every Yang--Baxter-type test passed. The relation that fails is the one that ties an R-matrix to a coproduct, $R\,\Delta(X)=\Delta^{\mathrm{op}}(X)\,R$. Listing~\ref{lst:coproduct} forms this residual for both matrices of Eq.~\eqref{eq:rmatrices} with public functions only, so that the same code can be pointed at any other $4\times4$ R-matrix.

\begin{lstlisting}[language=Python, caption={Testing an R-matrix against the coproduct: \code{scipost/examples/coproduct\_check.py}.}, label={lst:coproduct}]
import sympy as sp
from quantum_group import (build_representation, coproduct, swap_matrix,
                           R_matrix_V1, R_matrix_V1_coproduct, qybe_holds,
                           is_zero_matrix)

D = coproduct(build_representation(1))   # Delta(X) on V_1 (x) V_1, 4x4
P = swap_matrix(2)                       # ordinary flip of the two factors

def intertwining_residual(R, X):
    """R Delta(X) - Delta^op(X) R, with Delta^op(X) = P Delta(X) P."""
    return R * D[X] - P * D[X] * P * R

for name, R in [("R_matrix_V1", R_matrix_V1()),
                ("R_matrix_V1_coproduct", R_matrix_V1_coproduct())]:
    print(name, "- QYBE:", qybe_holds(R))
    for X in ("E", "F", "K"):
        res = intertwining_residual(R, X)
        nonzero = {ij: sp.factor(x) for ij, x in sorted(res.todok().items())
                   if sp.simplify(x) != 0}
        print(f"  {X}:", "zero" if is_zero_matrix(res) else nonzero)
\end{lstlisting}
Output (\code{scipost/\allowbreak examples/\allowbreak coproduct\_check\_expected.txt}; long lines wrapped here):
\begin{lstlisting}[style=output]
R_matrix_V1 - QYBE: True
  E: {(0, 1): (q - 1)*(q + 1), (0, 2): -(q - 1)*(q + 1)/q,
      (1, 3): (q - 1)*(q + 1)/q**2, (2, 3): -(q - 1)*(q + 1)/q}
  F: {(1, 0): (q - 1)*(q + 1)/q**2, (2, 0): -(q - 1)*(q + 1)/q,
      (3, 1): (q - 1)*(q + 1), (3, 2): -(q - 1)*(q + 1)/q}
  K: zero
R_matrix_V1_coproduct - QYBE: True
  E: zero
  F: zero
  K: zero
\end{lstlisting}
Both matrices pass the QYBE, but only $R_{21}$ intertwines $\Delta$. The residual of $R$ is informative in three ways. It vanishes for $K$, as it must for any weight-preserving matrix, so the problem is not the diagonal part. Every nonzero entry carries the factor $(q-1)(q+1)$, so all of them vanish at $q=\pm1$, where $q-q^{-1}=0$ and $R$ and $R_{21}$ coincide; the discrepancy is a genuine $q$-deformation effect. And every nonzero entry couples one of the middle basis vectors $v_0\otimes v_1$, $v_1\otimes v_0$ (indices 1 and 2), the block in which $R$ and $R_{21}$ differ, to $v_0\otimes v_0$ or $v_1\otimes v_1$. The same residual with $\Delta$ and $\Delta^{\mathrm{op}}$ exchanged vanishes for $R$, confirming that $R$ belongs to the opposite coproduct. The package now exposes this check as \code{intertwining\_residual\_V1}, and the test suite keeps the failing case as a negative control.

\subsection{Locating a transcription error}
\label{sec:ex-typo}

A matrix typed from a paper can be passed to \code{qybe\_residual} directly. Listing~\ref{lst:typo} introduces a single sign error in $R$ of Eq.~\eqref{eq:rmatrices}, writing $q+q^{-1}$ for $q-q^{-1}$.

\begin{lstlisting}[language=Python, caption={A mistyped entry: \code{scipost/examples/transcription\_check.py}.}, label={lst:typo}]
import sympy as sp
from quantum_group import q, qybe_residual, is_zero_matrix

# R-matrix typed from a source with one mistyped sign: q + 1/q at (1, 2),
# where the source (the package's R_matrix_V1) has q - 1/q.
R_typed = sp.Matrix([[q, 0, 0,       0],
                     [0, 1, q + 1/q, 0],
                     [0, 0, 1,       0],
                     [0, 0, 0,       q]])
Y = qybe_residual(R_typed)                 # exact 8x8 residual
print("QYBE holds:", is_zero_matrix(Y))
print("nonzero entries:", {ij: sp.factor(x) for ij, x in sorted(Y.todok().items())
                           if sp.simplify(x) != 0})
\end{lstlisting}
Output (\code{scipost/examples/transcription\_check\_expected.txt}):
\begin{lstlisting}[style=output]
QYBE holds: False
nonzero entries: {(1, 4): -2*(q**2 + 1)/q, (3, 6): 2*(q**2 + 1)/q}
\end{lstlisting}
Only two of the 64 entries are nonzero. In the basis $e_a\otimes e_b\otimes e_c\mapsto 4a+2b+c$ of $V_1^{\otimes3}$, they couple $v_0\otimes v_0\otimes v_1$ with $v_1\otimes v_0\otimes v_0$ and $v_0\otimes v_1\otimes v_1$ with $v_1\otimes v_1\otimes v_0$, so the search for the error can start from the entries that move a $v_1$ between the first and third factors. A floating-point check at a sample value of $q$ would report a failure as well, but not the exact form of the error.

\subsection{Clebsch--Gordan decomposition}
\label{sec:ex-cg}

Listing~\ref{lst:cg} decomposes $V_2\otimes V_1\cong V_3\oplus V_1$ by computing one highest-weight vector per summand in the tensor-product basis $v_i\otimes w_j\mapsto 2i+j$.

\noindent\begin{minipage}{\linewidth}
\begin{lstlisting}[language=Python, caption={Highest-weight vectors of $V_2\otimes V_1$: \code{scipost/examples/clebsch\_gordan.py}.}, label={lst:cg}]
from quantum_group import (build_representation, tensor_product,
                           find_highest_weight_vectors, cg_summands)

T = tensor_product(build_representation(2), build_representation(1))
print("summands V_k, k =", cg_summands(2, 1))
for k, v in find_highest_weight_vectors(T):
    print(f"k = {k}:", list(v))
\end{lstlisting}
\end{minipage}
Output (\code{scipost/examples/clebsch\_gordan\_expected.txt}):
\begin{lstlisting}[style=output]
summands V_k, k = [1, 3]
k = 3: [1, 0, 0, 0, 0, 0]
k = 1: [0, (-q**2 - 1)/q**3, 1, 0, 0, 0]
\end{lstlisting}
The highest-weight vector of the $V_1$ summand is $v_1\otimes w_0-q^{-3}(q^2+1)\,v_0\otimes w_1$; its $q$-dependent coefficient is fixed by the coproduct~\eqref{eq:coproduct}: with $\Delta^{\mathrm{op}}$ it would be $-(q^2+1)$ instead of $-q^{-3}(q^2+1)$, which is why the coproduct formulas are stated in the documentation of the \code{hopf} and \code{tensor} modules and in the README.

\subsection{Graded identities on four factors and a reordered basis}
\label{sec:ex-graded}

Listing~\ref{lst:graded4} runs two further graded checks for $\GL$. The first two lines verify the graded Yang--Baxter equation on every triple of factors of $V^{\otimes4}$, with all placements $R_{ij}$ built by graded adjacent swaps, and the commutation of the operators $R_{12}$ and $R_{34}$ on disjoint factors, an $81\times81$ identity. The second part tests the parity list as an argument. It rewrites the R-matrix~\eqref{eq:rgl21} in the reordered basis $f=(e_2,e_0,e_1)$, in which the odd vector comes first, and repeats the graded Yang--Baxter check once with the parities $(1,0,0)$ of the new basis and once with the stale list $(0,0,1)$ of the old one.

\begin{lstlisting}[language=Python, caption={Graded identities on $V^{\otimes4}$ and in a reordered basis: \code{scipost/examples/graded\_four\_factors.py}.}, label={lst:graded4}]
import sympy as sp
from quantum_group import (R_matrix_GLq21, embed_R_in_tensor_power,
                           local_ybe_on_four_tensor_GLq21,
                           braid_far_commutativity_residual_GLq21, is_zero_matrix)
from quantum_group.linalg import kron

print("local YBE on V^(x)4:", local_ybe_on_four_tensor_GLq21())
print("R12 R34 = R34 R12:", is_zero_matrix(braid_far_commutativity_residual_GLq21()))

# Reorder the basis so that the odd vector comes first: f = (e_2, e_0, e_1).
perm = [2, 0, 1]
Pi = sp.Matrix(3, 3, lambda k, i: int(perm[k] == i))     # new = Pi * old
R_new = kron(Pi, Pi) * R_matrix_GLq21() * kron(Pi, Pi).T
for parity in ([1, 0, 0], [0, 0, 1]):                     # updated / stale
    R12, R13, R23 = (embed_R_in_tensor_power(R_new, ij, 3, parity)
                     for ij in [(0, 1), (0, 2), (1, 2)])
    Y = R12 * R13 * R23 - R23 * R13 * R12
    print("parities", parity, "graded YBE holds:", is_zero_matrix(Y),
          "| nonzero entries:", sum(sp.simplify(x) != 0 for x in Y))
\end{lstlisting}
Output (\code{scipost/\allowbreak examples/\allowbreak graded\_four\_factors\_expected.txt}; the first line is wrapped here):
\begin{lstlisting}[style=output]
local YBE on V^(x)4: {(0, 1, 2): True, (0, 1, 3): True, (0, 2, 3): True, (1, 2, 3): True}
R12 R34 = R34 R12: True
parities [1, 0, 0] graded YBE holds: True | nonzero entries: 0
parities [0, 0, 1] graded YBE holds: False | nonzero entries: 4
\end{lstlisting}
With the parity list that matches the basis, the graded Yang--Baxter equation holds in the reordered basis as well; the Koszul signs are recomputed from the list, not from a fixed position of the odd vector. Reordering the basis but keeping the old parities, a mistake of the same kind as the transcription errors of Sec.~\ref{sec:intro}, leaves four nonzero entries in the $27\times27$ residual. The same call therefore serves any finite-dimensional super vector space whose parities are known.

\section{Testing and benchmarks}
\label{sec:tests}

\subsection{Test suite and continuous integration}

Every change is checked by one command, \code{python -m pytest}, which runs the test suite with pytest~\cite{pytest} together with the package doctests. At version 1.1.1, as at 1.1.0, the suite comprises 198 test cases in ten files and four doctests (202 in total, all passing). GitHub Actions runs it on every push and pull request on Python 3.10, 3.11, 3.12 and 3.13; once with the oldest supported dependency versions (SymPy 1.10.1, NetworkX 2.6.3, Matplotlib 3.5.3 on Python 3.10); and, after building and installing a wheel, runs the README quickstart and Listings~\ref{lst:example}--\ref{lst:crosscheck} from outside the source tree, comparing each listing's output with its stored copy. A further job repeats the wheel installation, the quickstart, Listing~\ref{lst:example} and the full test suite on Linux, macOS and Windows runners (Python 3.12), and the job with the oldest dependencies also runs Listings~\ref{lst:example}--\ref{lst:crosscheck}; their outputs are the same with SymPy 1.14 on Python 3.12 and with SymPy 1.10.1 on Python 3.10. Table~\ref{tab:qc} summarizes what the tests establish. Besides algebraic identities, the suite benchmarks the implementation against published reference values: every nonzero entry of the $\GL$ R-matrix is compared with the printed source~\cite{celik2021}, and the normalization of the $\Uq$ R-matrix is recomputed from Drinfeld's universal R-matrix~\cite{drinfeld1987}. Section~\ref{sec:crosscheck} adds a comparison with independently computed data from GAP's QuaGroup.

\begin{table}[t]
\centering
\caption{Evidence provided by the test suite. Entries in the last column abbreviate \code{tests/test\_\textit{name}.py}. ``Negative control'' marks tests that apply a known wrong convention and require a nonzero residual.}
\label{tab:qc}
\small
\begin{tabularx}{\textwidth}{@{}>{\raggedright\arraybackslash}p{0.17\textwidth}>{\raggedright\arraybackslash}X>{\raggedright\arraybackslash}p{0.235\textwidth}@{}}
\toprule
Evidence & What is checked & Test file \\
\midrule
Exact entries & Every nonzero entry of the $9\times9$ $\GL$ R-matrix against the display on printed page 261 of Ref.~\cite{celik2021}; exact entries of the $4\times4$ $\Uq$ R-matrix; shift and weight structure of $V_n$ & \texttt{audit\_regressions}, \texttt{r\_matrix}, \texttt{representations} \\
Coproduct convention & $R_{21}\Delta(X)=\Delta^{\mathrm{op}}(X)R_{21}$ for $X=E,F,K,K^{-1}$; negative control: the historical matrix fails for $E$ and $F$ & \texttt{audit\_regressions} \\
Clebsch--Gordan bases & Complete change of basis for $V_1\otimes V_1$, $V_2\otimes V_2$, $V_3\otimes V_2$: invertible over $\mathbb{Q}(q)$ and intertwining all four generators; braid eigenvalues $q$ and $-q^{-1}$ on the $V_2$ and $V_0$ summands & \texttt{audit\_regressions}, \texttt{tensor} \\
Graded embeddings & All six placements $R_{ij}$ on $V^{\otimes4}$ equal an independently written signed basis-action formula that builds no swap matrix; disjoint placements commute & \texttt{audit\_regressions}, \texttt{supergroup\_gl21} \\
Grading & Graded YBE residual ($27\times27$) is zero; negative control: with the ordinary flip exactly two entries are nonzero, with their values fixed & \texttt{supergroup\_gl21}, \texttt{audit\_regressions} \\
Exact parameters & No floating-point atoms for integer, rational, $q=-1$ and $q=i$ inputs; Laurent-polynomial $q$-integers and $q$-binomials agree with generic formulas after substitution; root-of-unity orders validated & \texttt{audit\_regressions}, \texttt{limits} \\
Algebraic identities & Defining relations on $V_n$; coassociativity, counit and antipode on generators; QYBE, braid and Hecke relations & \texttt{relations}, \texttt{hopf}, \texttt{r\_matrix} \\
Input validation & Every validation path raises the documented error; shared helpers agree with \code{sympy.kronecker\_product} & \texttt{validation\_and\_\allowbreak linalg} \\
\bottomrule
\end{tabularx}
\end{table}

\subsection{Negative controls and independent re-derivations}

A suite that only checks that correct identities hold cannot show that the checks are able to fail. Several tests therefore apply a known wrong convention (the opposite coproduct, the ordinary instead of the graded flip, a mutated counit) and require the residual to be nonzero, in some cases at specified entries with specified values: the two nonzero entries of the ungraded residual in Listing~\ref{lst:example} are fixed by a test, and the failure of the historical R-matrix to intertwine the coproduct (Listing~\ref{lst:coproduct}) is kept as a test requiring a nonzero residual. These tests fix the scientific content of the conventions, not only the behaviour of the code. Two further tests are independent re-derivations: the four-factor placements are compared with a formula that acts on basis vectors directly, and the normalization of the fundamental R-matrix is recomputed from Drinfeld's universal R-matrix~\cite{drinfeld1987}.

\subsection{Cross-check against GAP's QuaGroup}
\label{sec:crosscheck}

The tests above compare the package with published formulas. To compare it with existing software as well, the module matrices of $V_1$, $V_2$, $V_3$, the action of $\Uq$ on $V_n\otimes V_n$ ($n=1,2$) and the R-matrices \code{RMatrix}$(V_n)$ ($n=1,2$) were computed with QuaGroup~1.8.4~\cite{quagroup,degraaf2001} in GAP~4.12.1~\cite{gap}, the package behind SageMath's \code{QuantumGroup} interface~\cite{sagemath}, and exported by the GAP script \code{scipost/crosscheck/quagroup\_export.g}. The exported file (\code{quagroup\_data.txt}, in the repository) is read by Listing~\ref{lst:crosscheck}, which tests QuaGroup's data with the residual functions of \pkg{}. GAP is needed only to regenerate the data.

\begin{lstlisting}[language=Python, caption={Testing QuaGroup's data with \pkg: \code{scipost/crosscheck/quagroup\_crosscheck.py}.}, label={lst:crosscheck}]
import re
import sympy as sp
from quantum_group import (q, build_representation, coproduct, swap_matrix,
                           verify_on_representation, R_matrix_V1_coproduct,
                           qybe_holds, braid_relation_residual, is_zero_matrix)
from quantum_group.linalg import kron

G = {}  # matrices exported from QuaGroup by quagroup_export.g
for line in open(__file__.replace("quagroup_crosscheck.py", "quagroup_data.txt")):
    if m := re.match(r"(\w+) = (\[.*\])$", line.strip()):
        G[m[1]] = sp.Matrix(sp.sympify(m[2].replace("^", "**"), {"q": q}))
Z = is_zero_matrix

for n in (1, 2, 3):                      # QuaGroup's V_n vs the package's V_n
    E, F, K, Ki = (G[f"{X}_V{n}"] for X in ("E", "F", "K", "K_inv"))
    ok = all(c.holds for c in verify_on_representation(E, F, K, Ki).values())
    D = sp.diag(*[1 / sp.prod([sum(q**(k - 1 - 2*i) for i in range(k))
                               for k in range(1, j + 1)]) for j in range(n + 1)])
    r = build_representation(n)          # D = diag(1/[j]!): divided powers
    print(f"V_{n}: relations {ok}, equal after D:",
          Z(D * E * D.inv() - r.E) and Z(D * F * D.inv() - r.F) and K == r.K)

for n in (1, 2):                         # QuaGroup's V_n (x) V_n and R-matrix
    E, F, K, Ki = (G[f"{X}_V{n}"] for X in ("E", "F", "K", "K_inv"))
    I, P = sp.eye(n + 1), swap_matrix(n + 1)
    Delta = {"E": kron(E, I) + kron(K, E), "F": kron(F, Ki) + kron(I, F),
             "K": kron(K, K)}
    print(f"V_{n} (x) V_{n}: action = package coproduct:",
          all(Z(Delta[X] - G[f"{X}_T{n}"]) for X in "EFK"))
    for label, th in (("as printed", G[f"RMatrix_V{n}"]),
                      ("transposed", G[f"RMatrix_V{n}"].T)):
        commutes = all(Z(th * P * G[f"{X}_T{n}"] - G[f"{X}_T{n}"] * th * P)
                       for X in ("E", "F", "K"))
        print(f"  RMatrix {label}: QYBE {qybe_holds(th)},",
              f"braid {Z(braid_relation_residual(th * P))},",
              f"theta P commutes with Delta {commutes}")

th, R21 = G["RMatrix_V1"].T, R_matrix_V1_coproduct()
print("theta = q * R_21(1/q):", Z(th - q * R21.subs(q, 1 / q)))
D2 = sp.diag(1, 1, 1 / (q + 1 / q))
th2 = kron(D2, D2) * G["RMatrix_V2"].T * kron(D2, D2).inv()   # package basis
Rc2, cp = th2 * swap_matrix(3), coproduct(build_representation(2))
print("V_2 R-matrix in package basis: QYBE", qybe_holds(th2), "| commutes",
      all(Z(Rc2 * cp[X] - cp[X] * Rc2) for X in ("E", "F", "K")),
      "| eigenvalues", {sp.factor(k): v for k, v in Rc2.eigenvals().items()})
\end{lstlisting}
Output (\code{scipost/\allowbreak crosscheck/\allowbreak quagroup\_crosscheck\_expected.txt}; the last line is wrapped here):
\begin{lstlisting}[style=output]
V_1: relations True, equal after D: True
V_2: relations True, equal after D: True
V_3: relations True, equal after D: True
V_1 (x) V_1: action = package coproduct: True
  RMatrix as printed: QYBE True, braid True, theta P commutes with Delta False
  RMatrix transposed: QYBE True, braid True, theta P commutes with Delta True
V_2 (x) V_2: action = package coproduct: True
  RMatrix as printed: QYBE True, braid True, theta P commutes with Delta False
  RMatrix transposed: QYBE True, braid True, theta P commutes with Delta True
theta = q * R_21(1/q): True
V_2 R-matrix in package basis: QYBE True | commutes True
  | eigenvalues {q**(-2): 5, -q**2: 3, q**4: 1}
\end{lstlisting}
The two systems agree on the mathematics and differ in three conventions, each of which the residuals detect. (i)~QuaGroup's modules satisfy the relations~\eqref{eq:relations} and have the same $K$, but its basis is the divided-power basis $F^{(k)}v_0=F^kv_0/[k]!$; after conjugation with $D=\mathrm{diag}(1/[k]!)$ its $E$ and $F$ equal the package's. The action on $V_n\otimes V_n$ is exactly the coproduct~\eqref{eq:coproduct}. (ii)~GAP writes matrices in the row-vector convention; QuaGroup's source builds \code{RMatrix} row by row from the images of the basis vectors. Read as printed, \code{RMatrix}$(V_n)$ still satisfies the QYBE, which is invariant under transposition, and hence the braid relation, but $\theta P$ does not commute with the coproduct; only the intertwining test exposes the misreading. Transposed, the matrix passes all three tests. (iii)~For $V_1$ the transposed matrix is $\theta=q\,R_{21}(q^{-1})$, so that $\theta P=q\,\check R^{-1}$ with $\check R$ the braiding of \code{R\_check\_V1\_coproduct}: QuaGroup's braiding is the inverse of the package's, rescaled by $q$, and has eigenvalues $1$ (multiplicity three) and $-q^2$ instead of $q$ and $-q^{-1}$. Finally, QuaGroup's $9\times9$ R-matrix of $V_2$, which \pkg{} does not provide, satisfies the QYBE and commutes with the package coproduct once transported to the package basis; its braid eigenvalues $q^{-2}$, $-q^2$ and $q^4$ have multiplicities $5$, $3$ and $1$, the dimensions of the summands $V_4$, $V_2$ and $V_0$ of $V_2\otimes V_2$. External R-matrices can thus be brought into the package's conventions and verified there, which is the extension route described in Sec.~\ref{sec:scope}.

\subsection{Performance benchmarks}
\label{sec:benchmarks}

Operators on $V^{\otimes n}$ are stored as dense symbolic matrices, $3^n\times3^n$ for $\GL$, which keeps every entry inspectable but grows exponentially. Table~\ref{tab:bench} reproduces the recorded benchmark (\code{benchmarks/results.md}), generated by the command given in Sec.~\ref{sec:install}. Tensor powers $n\geq6$ have not been benchmarked. A sparse or tensor-network backend could reach larger $n$, at the cost of the entry-by-entry inspectability that motivates the design; the calculations the package targets involve at most four or five factors.

\begin{table}[t]
\centering
\caption{Wall-clock times (mean and minimum of three runs after a first, separately reported run) from \code{benchmarks/results.md}; the smallest case is dominated by timing noise (standard deviation 0.080\,s). Environment: Python 3.11, SymPy 1.14, one core of a 2.1\,GHz Intel Xeon processor, 27 September 2026.}
\label{tab:bench}
\small
\begin{tabularx}{\textwidth}{@{}>{\raggedright\arraybackslash}Xlrrr@{}}
\toprule
Computation & $n$ & Matrix size & Mean (s) & Min (s) \\
\midrule
QYBE of $R$ for $\Uq$ (\code{qybe\_holds}) & 3 & $8\times8$ & 0.075 & 0.028 \\
Graded YBE residual for $\GL$ & 3 & $27\times27$ & 0.025 & 0.024 \\
Graded YBE, residual and zero test & 3 & $27\times27$ & 0.060 & 0.055 \\
All placements $R_{ij}$ (\code{all\_Rij\_GLq21}) & 4 & $81\times81$ & 0.977 & 0.954 \\
Local YBE on every triple (\code{local\_ybe\_on\_four\_tensor\_GLq21}) & 4 & $81\times81$ & 1.556 & 1.540 \\
All placements $R_{ij}$ (\code{all\_Rij\_GLq21}) & 5 & $243\times243$ & 14.616 & 14.276 \\
Graded YBE on one triple, incl.\ embedding & 5 & $243\times243$ & 15.221 & 14.799 \\
\bottomrule
\end{tabularx}
\end{table}

\section{Scope, limitations and extensions}
\label{sec:scope}

\par\needspace{3\baselineskip}\smallskip\noindent\emph{Scope and limitations.} A zero residual establishes an identity on the particular finite-dimensional representation used; it is not a proof in the abstract Hopf algebra, and the package is not a theorem prover. It does not implement general quantum groups. For $\GL$ it implements the R-matrix, the graded Yang--Baxter equation and tensor placements, not the RTT relations, the Hopf superalgebra structure or the Gauss decomposition of Ref.~\cite{celik2021}. It verifies the Hecke relation but does not compute a Markov trace or the Jones polynomial. It works with constant R-matrices; spectral-parameter-dependent solutions of the Yang--Baxter equation are not built in. Its root-of-unity helpers substitute special values into the generic construction; they do not implement the small quantum group. The crystal graph $B(n)$ is a combinatorial model rather than a crystal-basis algorithm.

\par\needspace{3\baselineskip}\smallskip\noindent\emph{Reuse.} The extension points are ordinary Python functions that accept and return SymPy matrices, and the package is small enough to read in full. Any $d^2\times d^2$ matrix can be passed to \code{qybe\_residual} or \code{braid\_relation\_residual} (Sec.~\ref{sec:ex-typo}); the intertwining test of Listing~\ref{lst:coproduct} works for any R-matrix, module or coproduct formula; and \code{verify\_on\_representation} accepts any four square matrices, so alternative bases or normalizations of $V_n$ can be checked directly. \code{super\_permutation\_matrix} and \code{embed\_R\_in\_tensor\_power} take an arbitrary parity list and can place any even operator on factors $(i,j)$ of $V^{\otimes n}$ for other $\mathbb{Z}_2$-graded spaces, such as $\mathbb{C}^{m|n}$ with a different ordering of even and odd basis vectors. Because each check returns an exact matrix, an identity from a paper or thesis becomes a pytest test in a few lines, with a known wrong convention kept alongside it as a negative control; \code{tests/test\_audit\_regressions.py} contains worked examples of this pattern. Outputs can be exported with \code{sympy.latex} or compared against matrices computed elsewhere, for example in SageMath, by forming a residual.

\par\needspace{3\baselineskip}\smallskip\noindent\emph{Where reuse requires new work.} Moving from $\mathfrak{sl}_2$ to general semisimple Lie algebras would need root data, PBW bases and higher-rank modules, which QuaGroup and SageMath already provide. Implementing the RTT relations or the Hopf superalgebra of $\GL$ would need a noncommutative graded algebra rather than matrices. The dense matrix representation limits practical calculations to about five tensor factors of a three-dimensional space. These directions are possible but would be substantial new software rather than reuse of the existing layer. A lighter route is to import explicit matrices computed elsewhere, as done for the $V_2$ R-matrix of QuaGroup in Sec.~\ref{sec:crosscheck}, and verify them with the existing residual and embedding functions.

\section{Conclusion}
\label{sec:conclusion}

\pkg{} turns convention-sensitive identities of quantum-group theory into exact residual matrices that can be inspected entry by entry and kept as regression tests. Its value lies less in the range of algebras it covers, which is deliberately small, than in making the conventions explicit: basis order, coproduct and grading are fixed in documented functions, the Koszul sign rule is encoded once for every tensor placement, and every check can fail informatively. The thesis audit of Sec.~\ref{sec:ex-coproduct} shows that this matters in practice: an R-matrix that passes every Yang--Baxter-type test can still belong to the wrong coproduct, and only a residual tied to the coproduct reveals it. The comparison with QuaGroup in Sec.~\ref{sec:crosscheck} shows the same effect across software: the two systems agree on the mathematics, and the residuals locate exactly where their conventions for bases, matrix orientation and braiding differ. Natural next steps are further explicit rank-one and graded R-matrices, which reuse the existing residual and embedding layer unchanged, and a sparse backend for larger tensor powers.

\section*{Acknowledgments}

The author thanks Prof.\ Dr.\ Salih Çelik for supervising the undergraduate thesis from which this software originated.

\section*{Declarations}

\noindent\textbf{Author contributions:} T.B.T. is the sole author: T.B.T. designed and wrote the original implementation; decided on, reviewed and integrated all later changes to the software, including AI-assisted ones; and prepared and approved the manuscript.

\medskip
\noindent\textbf{Funding:} This work received no specific funding.

\medskip
\noindent\textbf{Competing interests:} The author has no competing interests to declare.

\medskip
\noindent\textbf{Code and data availability:} The software, the scripts of all listings with their expected outputs and the exported QuaGroup data are in the repository and its release archive; see the paragraph \emph{Code availability} in Section~\ref{sec:install}.

\medskip
\noindent\textbf{Use of generative AI:} No generative AI tools were used during the original development and thesis work (April--June 2026). In September 2026, Claude Code (Anthropic) was used for packaging and continuous-integration configuration, the benchmark script, the shared matrix helpers, input validation and their tests, the sample-verification example, the worked-example and cross-check scripts of this paper, documentation, translating Turkish docstrings and comments into English, a mathematical and technical audit that led to the parameter-handling corrections, the coproduct-compatible R-matrix functions and the audit regression tests, and drafting of this manuscript. The model was Claude Opus 5.5, except for part of the manuscript and submission-material preparation, which used Claude Sonnet 5; commits produced with Claude Code are authored or co-authored by ``Claude'' in the Git history, with the model named in the commit trailer. No other AI tool was used. Every AI-assisted code change is exercised by the test suite, and every AI-assisted mathematical claim was checked by an executable symbolic test or against the primary source, in particular the published Çelik--Çelik article. The author reviewed and validated all AI-assisted output, made the scientific and design decisions, and takes responsibility for the software and this paper.

\appendix

\section{Reference of the public interface}
\label{app:api}
\begingroup\raggedright

This appendix lists the main public functions of \pkg{} 1.1.1 with their arguments and return values; the docstrings give the full conventions. Unless stated otherwise, functions are imported from the package namespace (\texttt{from quantum\_group import ...}), matrices are exact \texttt{sympy.Matrix} objects in the column convention, and tensor-product bases are ordered lexicographically with the left factor slowest. Functions whose result depends on $q$ take a final optional argument \texttt{q\_sym}, which defaults to the package symbol \texttt{q} and is omitted below. (\texttt{q} is declared nonzero; SymPy's \texttt{nonzero} assumption also makes it real, which does not affect identities between rational functions of $q$.) A \texttt{*\_holds} function returns \texttt{True} exactly when the corresponding \texttt{*\_residual} matrix passes \texttt{is\_zero\_matrix}.

\par\needspace{3\baselineskip}\smallskip\noindent\emph{Records.}
\begin{itemize}\setlength{\itemsep}{1pt}\setlength{\parskip}{0pt}
\item \texttt{Representation}: fields \texttt{n}, \texttt{dim}, \texttt{E}, \texttt{F}, \texttt{K}, \texttt{K\_inv}, \texttt{weights}.
\item \texttt{TensorRepresentation}: fields \texttt{m}, \texttt{n}, \texttt{dim}, \texttt{E}, \texttt{F}, \texttt{K}, \texttt{K\_inv} (action on $V_m\otimes V_n$).
\item \texttt{RelationCheck}: fields \texttt{name}, \texttt{holds}, \texttt{residual}; \texttt{HopfAxiomCheck}: fields \texttt{name}, \texttt{holds}, \texttt{detail}.
\end{itemize}

\par\needspace{3\baselineskip}\smallskip\noindent\emph{Modules and relations (Sec.~\ref{sec:uq}).}
\begin{itemize}\setlength{\itemsep}{1pt}\setlength{\parskip}{0pt}
\item \texttt{build\_representation(n)} $\to$ \texttt{Representation}: the module $V_n$ of Eq.~\eqref{eq:module}; \texttt{n} a nonnegative integer.
\item \texttt{verify\_on\_representation(E, F, K, K\_inv)} $\to$ \texttt{dict}: keys \texttt{"R1"}--\texttt{"R4"} for the relations~\eqref{eq:relations}, values \texttt{RelationCheck}. Accepts any four square matrices of equal size; rejects $q\in\{0,1,-1\}$.
\item \texttt{all\_relations\_hold(checks)} $\to$ \texttt{bool}; \texttt{symbolic\_relations()} $\to$ list of \texttt{(name, lhs, rhs)} for display.
\end{itemize}

\par\needspace{3\baselineskip}\smallskip\noindent\emph{Hopf structure (Eqs.~\eqref{eq:coproduct}, \eqref{eq:hopf}).}
\begin{itemize}\setlength{\itemsep}{1pt}\setlength{\parskip}{0pt}
\item \texttt{coproduct(rep)} $\to$ \texttt{dict} of the matrices $\Delta(X)$ on $V\otimes V$, keys \texttt{"E"}, \texttt{"F"}, \texttt{"K"}, \texttt{"K\_inv"}.
\item \texttt{counit()} $\to$ \texttt{dict} of $\varepsilon(X)$; \texttt{antipode(rep)} $\to$ \texttt{dict} of the matrices $S(X)$ on $V$.
\item \texttt{verify\_coassociativity(rep)} and \texttt{verify\_antipode(rep)} $\to$ \texttt{dict} of \texttt{HopfAxiomCheck}.
\item \texttt{verify\_all\_hopf\_axioms(rep)} $\to$ the same checks, grouped by axiom.
\end{itemize}

\par\needspace{3\baselineskip}\smallskip\noindent\emph{Tensor products.}
\begin{itemize}\setlength{\itemsep}{1pt}\setlength{\parskip}{0pt}
\item \texttt{tensor\_product(repA, repB)} $\to$ \texttt{TensorRepresentation}, the action through $\Delta$.
\item \texttt{find\_highest\_weight\_vectors(T)} $\to$ list of \texttt{(k, vector)}, one highest-weight vector per summand $V_k$.
\item \texttt{cg\_summands(m, n)} $\to$ list of the $k$ in $V_m\otimes V_n\cong\bigoplus V_k$; \texttt{cg\_decomposition\_summary(m, n)} $\to$ \texttt{str}.
\item \texttt{swap\_matrix(d)} $\to$ the flip $P$ on $\mathbb{C}^d\otimes\mathbb{C}^d$; \texttt{kron\_list(mats)} $\to$ left-to-right Kronecker product; \texttt{kron(A, B)} (from \texttt{quantum\_group.linalg}).
\end{itemize}

\par\needspace{3\baselineskip}\smallskip\noindent\emph{\texorpdfstring{$\Uq$}{Uq(sl2)} R-matrices (Sec.~\ref{sec:conventions}).}
\begin{itemize}\setlength{\itemsep}{1pt}\setlength{\parskip}{0pt}
\item \texttt{R\_matrix\_V1()} $\to$ $R$ of Eq.~\eqref{eq:rmatrices} (intertwines $\Delta^{\mathrm{op}}$); \texttt{R\_check\_V1()} $\to$ $PR$.
\item \texttt{R\_matrix\_V1\_coproduct()} $\to$ $R_{21}=PRP$ (intertwines $\Delta$); \texttt{R\_check\_V1\_coproduct()} $\to$ $PR_{21}$.
\item \texttt{R\_check\_eigenvalues()} $\to$ \texttt{dict} eigenvalue $\to$ multiplicity; \texttt{hecke\_skein\_relation\_check()} $\to$ \texttt{dict} with the Boolean of Eq.~\eqref{eq:hecke}.
\item \texttt{intertwining\_residual\_V1(X)} $\to$ $R_{21}\Delta(X)-\Delta^{\mathrm{op}}(X)R_{21}$ for \texttt{X} in \texttt{"E"}, \texttt{"F"}, \texttt{"K"}, \texttt{"K\_inv"}; \texttt{intertwining\_holds\_V1(X)} $\to$ \texttt{bool}.
\end{itemize}

\par\needspace{3\baselineskip}\smallskip\noindent\emph{Generic residuals (any $d^2\times d^2$ matrix).}
\begin{itemize}\setlength{\itemsep}{1pt}\setlength{\parskip}{0pt}
\item \texttt{qybe\_residual(R)} $\to$ $R_{12}R_{13}R_{23}-R_{23}R_{13}R_{12}$ on $V^{\otimes3}$ (ungraded placement); \texttt{qybe\_holds(R)} $\to$ \texttt{bool}.
\item \texttt{braid\_relation\_residual(Rc)} $\to$ $\check R_{12}\check R_{23}\check R_{12}-\check R_{23}\check R_{12}\check R_{23}$ for a braided matrix \texttt{Rc}; \texttt{braid\_relation\_holds(Rc)} $\to$ \texttt{bool}.
\item \texttt{is\_zero\_matrix(M)} $\to$ \texttt{True} if every entry simplifies to exactly 0 (Sec.~\ref{sec:residual}).
\end{itemize}

\par\needspace{3\baselineskip}\smallskip\noindent\emph{Grading and \texorpdfstring{$\GL$}{GLq(2|1)} (Sec.~\ref{sec:grading}).}
\begin{itemize}\setlength{\itemsep}{1pt}\setlength{\parskip}{0pt}
\item \texttt{super\_parity\_gl21()} $\to$ \texttt{[0, 0, 1]}; \texttt{super\_permutation\_matrix(parity)} $\to$ $P_s$ of Eq.~\eqref{eq:superperm} for any parity list.
\item \texttt{embed\_R\_in\_tensor\_power(R, (i, j), n, parity)} $\to$ the placement $R_{ij}$ on $V^{\otimes n}$ (zero-based, $i<j$) by graded adjacent swaps; rejects odd operators, invalid positions and parity lists.
\item \texttt{R\_matrix\_GLq21()} $\to$ Eq.~\eqref{eq:rgl21}; \texttt{all\_Rij\_GLq21(n)} $\to$ \texttt{dict} $(i,j)\to R_{ij}$ on $V^{\otimes n}$.
\item \texttt{graded\_yang\_baxter\_residual\_GLq21()} $\to$ $27\times27$ residual; \texttt{graded\_yang\_baxter\_holds\_GLq21()} $\to$ \texttt{bool}.
\item \texttt{local\_ybe\_on\_four\_tensor\_GLq21()} $\to$ \texttt{dict} triple $\to$ \texttt{bool}; \texttt{braid\_far\_commutativity\_residual\_GLq21()} $\to$ $R_{12}R_{34}-R_{34}R_{12}$ ($81\times81$).
\end{itemize}

\par\needspace{3\baselineskip}\smallskip\noindent\emph{Specializations, crystals and plots.}
\begin{itemize}\setlength{\itemsep}{1pt}\setlength{\parskip}{0pt}
\item \texttt{q\_integer(n)}, \texttt{q\_factorial(n)}, \texttt{q\_binomial(n, k)} $\to$ Laurent polynomials in $q$; \texttt{classical\_limit(expr)} $\to$ the limit $q\to1$.
\item \texttt{root\_of\_unity\_substitution(rep, N)} $\to$ \texttt{dict} with $E^N$, $F^N$, $K^{2N}$ at $q=e^{2\pi i/N}$ ($N>2$), substituted into the generic matrices.
\item \texttt{build\_crystal(n)} $\to$ NetworkX \texttt{DiGraph} of $B(n)$ with string node labels; \texttt{crystal\_nodes(n)} $\to$ list of \texttt{CrystalNode} records (\texttt{index}, \texttt{weight}, \texttt{label}).
\item \texttt{e\_tilde(node, n)}, \texttt{f\_tilde(node, n)} $\to$ \texttt{CrystalNode} or \texttt{None}, for a \texttt{CrystalNode} from \texttt{crystal\_nodes(n)}.
\item \texttt{plot\_weight\_diagram(n)}, \texttt{plot\_crystal\_graph(n)} $\to$ Matplotlib \texttt{Figure}.
\end{itemize}
\par\endgroup

\bibliographystyle{amsplain}
\bibliography{references}

\end{document}